\documentclass[12]{amsart}
\usepackage{amsmath,amssymb,amsthm,color,enumerate,comment,centernot,enumitem,url,cite}
\usepackage{graphicx,relsize,bm}
\usepackage{mathtools}
\usepackage{array}

\makeatletter
\newcommand{\tpmod}[1]{{\@displayfalse\pmod{#1}}}
\makeatother

\newtheorem{thm}{Theorem}[section]

\theoremstyle{remark}

\theoremstyle{definition}

\theoremstyle{THM}

\def\ds{\displaystyle}

\def\FF {{\mathcal F}}

\def\Z {{\mathbb Z}}

\def\Q {{\mathbb Q}}

\def\Z {{\mathbb Z}}
\def\Q {{\mathbb Q}}

\def\Gal{{{\rm Gal} }}

\makeatletter
\@namedef{subjclassname@2020}{%
  \textup{2020} Mathematics Subject Classification}
\makeatother

\def\red#1 {\textcolor{red}{#1 }}
\def\blue#1 {\textcolor{blue}{#1 }}

\numberwithin{equation}{section}
\def\ds{\displaystyle}
\def\Z {{\mathbb Z}}

\newcommand{\Mod}[1]{\ (\mathrm{mod}\enspace #1)}
\newcommand{\mmod}[1]{\ \mathrm{mod}\enspace #1}
\begin{document}

\title[Monogenic Cyclic Quartic Trinomials]{Monogenic Cyclic Quartic Trinomials}


\author{Lenny Jones}
\address{Professor Emeritus, Department of Mathematics, Shippensburg University, Shippensburg, Pennsylvania 17257, USA}
\email[Lenny~Jones]{doctorlennyjones@gmail.com}

\date{\today}

\begin{abstract}
A monic polynomial $f(x)\in \Z[x]$ of degree $N$ is called \emph{monogenic} if $f(x)$ is irreducible over $\Q$ and $\{1,\theta,\theta^2,\ldots ,\theta^{N-1}\}$ is a basis for the ring of integers of $\Q(\theta)$, where $f(\theta)=0$. In this brief note, we prove that there exist exactly three distinct monogenic trinomials of the form $x^4+bx^2+d$ whose Galois group is the cyclic group of order 4.
\end{abstract}


\maketitle
\section{Introduction}\label{Section:Intro}

  We say that a monic polynomial $f(x)\in \Z[x]$ is \emph{monogenic} if $f(x)$ is irreducible over $\Q$ and $\{1,\theta,\theta^2,\ldots ,\theta^{\deg{f}-1}\}$ is a basis for the ring of integers $\Z_K$ of $K=\Q(\theta)$, where $f(\theta)=0$. When $f(x)$ is irreducible over $\Q$, we have \cite{Cohen}
\begin{equation} \label{Eq:Dis-Dis}
\Delta(f)=\left[\Z_K:\Z[\theta]\right]^2\Delta(K),
\end{equation}
where $\Delta(f)$ and $\Delta(K)$ denote the discriminants over $\Q$, respectively, of $f(x)$ and the number field $K$.
Thus, for irreducible $f(x)$, we have that $f(x)$ is monogenic if and only if
  $\Delta(f)=\Delta(K)$. We also say that any number field $K$ is \emph{monogenic} if there exists a power basis for $\Z_K$. We caution the reader that, while the monogenicity of $f(x)$ implies the monogenicity of $K=\Q(\theta)$, where $f(\theta)=0$,
 the converse is not necessarily true. A simple example is $f(x)=x^2-5$ and $K=\Q(\theta)$, where $\theta=\sqrt{5}$. Then, $\Delta(f)=20$ and $\Delta(K)=5$. Thus, $f(x)$ is not monogenic, but nevertheless, $K$ is monogenic since $\{1,(\theta+1)/2\}$ is a power basis for $\Z_K$. Observe then that $g(x)=x^2-x-1$, the minimal polynomial for $(\theta+1)/2$ over $\Q$, is monogenic.

This note was motivated by a recent question of Tristan Phillips (private communication) asking if it is possible to determine all distinct monogenic quartic trinomials that have Galois group $C_4$, the cyclic group of order 4. We consider two monogenic $C_4$-quartic trinomials $f(x)$ and $g(x)$ to be \emph{distinct} if $\Q(\alpha)\not \simeq\Q(\theta)$, where $f(\alpha)=0=g(\theta)$. In this note, we provide a partial answer to the question of Phillips by proving the following theorem.
\begin{thm}\label{Thm:Main}
The three trinomials 
\[x^4-4x^2+2,\quad x^4+4x^2+2\quad \mbox{and} \quad x^4-5x^2+5,\] are the only distinct 
trinomials of the form $f(x)=x^4+bx^2+d\in \Z[x]$ with $\Gal(f)\simeq C_4$.     
\end{thm}

\section{Preliminaries}\label{Section:Prelim}
The following theorem follows from results due to Kappe and Warren.
 
 \begin{thm}{\rm \cite{KW}}\label{Thm:KW1}
Let $f(x)=x^4+bx^2+d\in \Z[x]$. Then $f(x)$ is irreducible over $\Q$ with $\Gal(f)\simeq C_4$ if and only if 
\begin{equation}\label{Eq:cons}
d \mbox{ and } b^2-4d \mbox{ are not squares in $\Z$, but $d(b^2-4d)$ is a square in $\Z$.}
\end{equation}   
 \end{thm}
 
The next result is the specific case for our quartic situation of a ``streamlined" version of Dedekind's index criterion for trinomials that is due to Jakhar, Khanduja and Sangwan. We have used Swan's formula \cite{Swan} for the discriminant of an arbitrary trinomial $f(x)$ to calculate  $\Delta(f)$.
\begin{thm}{\rm \cite{JKS2}}\label{Thm:JKS}
Let $K=\Q(\theta)$ be an algebraic number field with $\theta\in \Z_K$, the ring of integers of $K$, having minimal polynomial $f(x)=x^{4}+bx^2+d$ over $\Q$. A prime factor $q$ of $\Delta(f)=2^{4}d(b^2-4d)^2$ does not divide $\left[\Z_K:\Z[\theta]\right]$ if and only if $q$ satisfies one of the following conditions:
\begin{enumerate}[font=\normalfont]
  \item \label{JKS:I1} when $q\mid b$ and $q\mid d$, then $q^2\nmid d$;
  \item \label{JKS:I2} when $q\mid b$ and $q\nmid d$, then
  \[\mbox{either } \quad q\mid b_2 \mbox{ and } q\nmid d_1 \quad \mbox{ or } \quad q\nmid b_2\left(-db_2^2-d_1^2\right),\]
  where $b_2=b/q$ and $d_1=\frac{d+(-d)^{q^j}}{q}$ with $q^j\mid\mid 4$;
  \item \label{JKS:I3} when $q\nmid b$ and $q\mid d$, then
  \[\mbox{either } \quad q\mid b_1 \mbox{ and } q\nmid d_2 \quad \mbox{ or } \quad q\nmid b_1d_2\left(-bb_1+d_2\right),\]
  where $b_1=\frac{b+(-b)^{q^e}}{q}$ with $q^e\mid\mid 2$, and $d_2=d/q$;
  \item \label{JKS:I4} when $q=2$ and $2\nmid bd$, then the polynomials
   \begin{equation*}
    H_1(x):=x^2+bx+d \quad \mbox{and}\quad H_2(x):=\dfrac{bx^2+d+\left(-bx-d\right)^2}{2}
   \end{equation*}
   are coprime modulo $2$;
         \item \label{JKS:I5} when $q\nmid 2bd$, then $q^2\nmid \left(b^2-4d\right)$.
   \end{enumerate}
\end{thm}

\section{The Proof of Theorem \ref{Thm:Main}}\label{Section:MainProof}

\begin{proof}
 Following Theorem \ref{Thm:KW1}, we assume conditions \eqref{Eq:cons} so that $f(x)$ is irreducible over $\Q$ with $\Gal(f)\simeq C_4$. Observe  that if $d<0$, then $d(b^2-4d)<0$, which contradicts the fact that $d(b^2-4d)$ is a square. Hence, $d>0$ and $b^2-4d>0$. Furthermore, since $d$ and $b^2-4d$ are not squares, but $d(b^2-4d)$ is a square, we deduce that $d\ge 2$ and $b^2-4d\ge 2$.
 
 We use Theorem \ref{Thm:JKS} to ``force" the monogenicity of $f(x)$. Let $q$ be a prime divisor of $d$. If $q\nmid (b^2-4d)$, then $q\nmid b$, and  $q^2\mid d$ since $d(b^2-4d)$ is a square. But then condition \eqref{JKS:I3} of Theorem \ref{Thm:JKS} is not satisfied since $q\mid d_2$. Therefore, $q\mid (b^2-4d)$, and so $q\mid b$. Note then that if $q^2\mid d$, then condition \eqref{JKS:I1} is not satisfied. Hence, $q\mid \mid d$ and therefore, $d$ is squarefree, $d\mid (b^2-4d)$ and $d\mid b$. 

 Suppose next that $q$ is a prime divisor of $b^2-4d$, such that $q\nmid d$. If $q\mid b$, then $q=2$ and 
 \begin{equation}\label{Eq:A}
 A:=d(b^2-4d)/4 \quad \mbox{is a square in $\Z$.}
 \end{equation} We examine condition \eqref{JKS:I2} of Theorem \ref{Thm:JKS} to see that
 \[d_1=\dfrac{d+(-d)^4}{2}\equiv \left\{\begin{array}{cl}
   1 \pmod{4} & \mbox{if $d\equiv 1 \pmod{4}$}\\[.5em]
   2 \pmod{4} & \mbox{if $d\equiv 3 \pmod{4}$.}
 \end{array}\right.\]
 Thus, the first statement under condition \eqref{JKS:I2} is satisfied if and only if 
 \begin{equation}\label{S1}
 (b \mmod{4},\ d\mmod{4})=(0,1),
 \end{equation}
  while the second statement under condition \eqref{JKS:I2} is satisfied if and only if
  \begin{equation}\label{S2}
  (b \mmod{4},\ d\mmod{4})=(2,3).
  \end{equation}
 In scenario \eqref{S1}, we have that $A\equiv 3 \pmod{4}$, while in scenario \eqref{S2}, we have that $A\equiv 2 \pmod{4}$, contradicting \eqref{Eq:A} in each scenario. Hence, $q\nmid b$ and $q\ge 3$. Since $q\nmid d$ and $d(b^2-4d)$ is a square, we must have that $q^2\mid (b^2-4d)$. But then condition \eqref{JKS:I5} of Theorem \ref{Thm:JKS} is not satisfied. Therefore, every prime divisor of $b^2-4d$ divides $d$. 
 
 Thus, to emphasize, we now have that $d$ is squarefree, and that $d$ and $b^2-4d$ have exactly the same prime divisors $p_1<p_2<\cdots <p_k$. 
  Hence, since $d(b^2-4d)$ is a square, we can write
 \begin{equation}\label{Eq:factors1}
 d(b^2-4d)=\left(\prod_{i=1}^kp_i\right)\left(b^2-4\left(\prod_{i=1}^kp_i\right)\right)=\prod_{i=1}^kp_i^{2e_i},
 \end{equation}
  for some integers $e_i\ge 1$. Then, from \eqref{Eq:factors1}, we have  
  \begin{equation*}\label{Eq:factors2}
    b^2=\left(\prod_{i=1}^kp_i\right)\left(\left(\prod_{i=1}^kp_i^{2e_i-2}\right)+4\right),
  \end{equation*}
 which implies that 
 \begin{equation}\label{Eq:divisibility}
 \prod_{i=1}^kp_i \quad \mbox{divides} \quad \left(\prod_{i=1}^kp_i^{2e_i-2}\right)+4.
 \end{equation} We see from \eqref{Eq:divisibility} that if some $e_i>1$, then $p_i\mid 4$ so that $i=1$ and $p_1=2$. In this case, we get that 
 \begin{equation}\label{Eq:b2}
 b^2=2\left(\prod_{i=2}^kp_i\right)\left(4^{e_1-1}+4\right)=\left\{
 \begin{array}{cl}
  2^4\ds \prod_{i=2}^kp_i & \mbox{if $e_1=2$}\\[.5em]
  2^3\left(\ds \prod_{i=2}^kp_i\right)\left(4^{e_1-2}+1\right) & \mbox{if $e_1\ge 3$.}
   \end{array}\right.
   \end{equation} The second case of \eqref{Eq:b2} is impossible since $b^2/8\equiv 1 \pmod{2}$. The first case of \eqref{Eq:b2} is viable provided $k=1$, so that $b^2=16$ and $d=2$. We then get the two trinomials 
   \[x^4-4x^2+2 \quad \mbox{and} \quad x^4+4x^2+2,\] which are both easily confirmed to be monogenic using Theorem \ref{Thm:JKS}.

The remaining possibility in \eqref{Eq:divisibility} when $e_i=1$ for all $i$ yields $k=1$ and $p_1=5$, so that $b^2=25$ and $d=5$. The two resulting trinomials are then
\[x^4+5x^2+5 \quad \mbox{and}\quad x^4-5x+5.\] Again, using Theorem \ref{Thm:JKS}, it is straightforward to verify that $x^4+5x^2+5$ is not monogenic (condition \eqref{JKS:I4} fails), while $x^4-5x^2+5$ is monogenic.

Thus, we have found exactly three monogenic cyclic trinomials 
\[x^4-4x^2+2,\quad x^4+4x^2+2\quad \mbox{and} \quad x^4-5x^2+5.\]
Note that
\begin{equation}\label{Eq:DisTri}
\Delta(x^4-4x^2+2)=\Delta(x^4+4x^2+2)=2^{11} \quad  \mbox{and} \quad \Delta(x^5-5x^2+5)=2^45^3.
\end{equation}
If any two of these three trinomials generate the same quartic field, then their discriminants must be equal since they are monogenic. Hence, we see immediately from \eqref{Eq:DisTri} that the quartic field generated by $x^5-5x^2+5$ is distinct from the other two quartic fields. However, equality of two discriminants is not sufficient to conclude that those trinomials generate isomorphic quartic fields.
 Indeed, since the field generated by $x^4-4x^2+2$ is real, while the field generated by $x^4+4x^2+2$ is non-real, we deduce that these two fields are in fact distinct. Alternatively, we can verify that these two fields are not isomorphic using MAGMA.
\end{proof}





\end{document}